%
%
%
%
%
%
%
\input amstex
\documentstyle{amsppt} 
\def\endth{\par\ifdim\lastskip<\bigskipamount
\removelastskip\penalty55
\bigskip\fi}

\def\ifundefined#1{\expandafter\ifx
\csname#1\endcsname\relax}
\newif\ifdevelop \developfalse

\newtoks\chnumber
\newtoks\sectionnumber
\newcount\equationnumber
\newcount\thnumber

\def\assignnumber#1#2{%
\ifundefined{#1}\relax\else\message{#1 already defined}\fi
\expandafter\xdef\csname#1\endcsname
{\if-\the\sectionnumber\else\the\sectionnumber.\fi\the#2}%
}%

\def\beginsektion #1 #2 {\vskip0pt plus.1\vsize\penalty-250
\vskip0pt plus-.1\vsize
\bigbreak\bigskip
\sectionnumber{#1} \equationnumber0\thnumber0
\noindent{\bf #1. #2}\par
\nobreak\medskip\noindent}
\def\NoBlackBoxes{\overfullrule0pt}

\def\nologo{\expandafter\let
\csname logo\string @\endcsname=\empty}
\def\:{:\allowbreak }
\def\eq#1{\relax
\global\advance\equationnumber by 1
\assignnumber{EN#1}\equationnumber
{\tenrm (\csname EN#1\endcsname)}
}

\def\eqtag#1{\relax\ifundefined{EN#1}
\message{EN#1 undefined}{\sl (#1)}%
\else(\csname EN#1\endcsname)\fi%
}

\def\thname#1{\relax
\global\advance\thnumber by 1
\assignnumber{TH#1}\thnumber
\csname TH#1\endcsname
}

\def\beginth#1 #2 {\bigbreak\noindent{\bf #1\enspace \thname{#2}}
    ---\hskip4pt}
\def\thtag#1{\relax\ifundefined{TH#1}
\message{TH#1 undefined}{\sl #1}%
  \else[\csname TH#1\endcsname]\fi}

\def\de{\delta}

\def\si{\sigma}

\def\De{\Delta}

\def\CC{{\Bbb C}}

\def\RR{{\Bbb R}}
\def\TT{{\Bbb T}}
\def\ZZ{{\Bbb Z}}
\def\FSA{\Cal A}
\def\FSB{\Cal B}
\def\FSC{\Cal C}

\def\FSH{\Cal H}

\def\FSU{\Cal U}

\def\FSZ{\Cal Z}

\def\gog{{\goth g}}
\def\gok{{\goth k}}
\def\gol{{\goth l}}

\def\goS{{\goth S}}
\def\id{\operatorname{id}}
\def\phi{\varphi}
\def\Aq{\FSA_q}
\def\Uq{\FSU_q}

\def\detq{\det\nolimits_q}
\def\Lp{L^+}
\def\Lm{L^-}
\def\Lpm{L^\pm}
\def\Hst{\FSH^{(\si,\tau)}}

\def\AT{\FSA(\TT)}

\def\End{\operatorname{End}}

\def\dd{\operatorname{d}}

\let\eps=\varepsilon

\let\ten=\otimes
\def\AW{AW}
\def\DKqsf{DK1}
\def\DKcqg{DK2}
\def\DN{DN}
\def\DZ{Dz}

\def\DR{Dr}
\def\JB{J}
\def\KWqsf{KW1}
\def\KWmaw{KW2}
\def\KS{KS}
\def\KU{KU}
\def\MCorth{M}
\def\NOmac{N}
\def\NM{NM}
\def\NS{NS}
\def\NYM{NYM}
\def\PO{P}
\def\RTF{RTF}
\def\RS{Ro}
\def\STjac{ST1}
\def\STaw{ST2}
\def\SK{SK}

\def\KV{KV}
\def\WZ{Wz}
\noindent
{\eightpoint \sl
To appear in:  Proceedings of a Workshop ``Special Functions, q-Series
 and Related Topics'', Toronto, June 19-23, 1995, Fields Inst. Comm.
}
\par\vskip1cm
\topmatter
\title
Multivariable Askey-Wilson polynomials\\
and quantum complex Grassmannians
\endtitle
\author 
Masatoshi Noumi${}^1$, Mathijs S.\ Dijkhuizen${}^1$, 
Tetsuya Sugitani${}^2$\endauthor
\affil ${}^1$ Department of Mathematics, Kobe University\\
Rokko, Kobe 657, Japan \\
${}^2$ Department of Mathematical Sciences, University of Tokyo\\
Komaba, Meguro-Ku, Tokyo 153, Japan \endaffil
\email msdz$\@$math.s.kobe-u.ac.jp,  
noumi$\@$math.s.kobe-u.ac.jp, sugitani$\@$ms.u-tokyo.ac.jp\endemail
\thanks The 
second author (M.S.\,D.)
acknowledges financial 
support by the Japan Society for the Promotion
of Science (JSPS) and the Netherlands Organization 
for Scientific Research (NWO). \endthanks
\abstract We present a one-parameter family of constant 
solutions of the reflection equation and define a family of 
quantum complex Grassmannians endowed with a 
transitive action of the quantum unitary group.
By computing the radial part of a suitable
Casimir operator, we identify the zonal spherical 
functions (i.e.\ infinitesimally bi-invariant 
matrix coefficients of finite-dimensional 
irreducible representations)
as multivariable Askey-Wilson polynomials containing two 
continuous and two discrete parameters. 
\endabstract
\date 15 September 1995 \enddate
\rightheadtext{Quantum complex Grassmannians}
\leftheadtext{M.\,Noumi, M.S.\,Dijkhuizen and T.\,Sugitani}
\endtopmatter
\NoBlackBoxes
\document
\beginsektion 0 {Introduction}
Every classical compact symmetric
space $M$ can be realized as an orbit of square matrices, the
action of the symmetry group $G$ being given 
by $X\mapsto TXT^t$ or
$X\mapsto TXT^\ast$ ($X\in M$, $T\in G$). 
This approach 
has proved to be fruitful in the study of quantum symmetric 
spaces as well. 
The proper
way to quantize the space of square matrices endowed
with either of
the above-mentioned group actions would be to
make use of commutation relations given by one of two types of 
so-called
{\sl reflection equations}. These equations pop up at various
places in mathematics, notably in the theory of quantum integrable 
systems (cf.\ \cite{\KU}). 
The quantized symmetric
space $M_q$ is then realized as the quantum orbit containing
a certain ``classical''  or $\CC$-valued point  in the quantum space of 
square matrices (i.e.\ a 
$\CC$-valued evaluation homomorphism on the quantized algebra 
of functions). 
These ``classical'' points are naturally 
in 1-1 correspondence with constant solutions $J$ of the
reflection equation. 

There is another, slightly different, 
but closely related, way
to explain the role of the reflection equation in the theory of quantum
symmetric spaces. 
Given any square matrix $J$ with complex coefficients,
there is a straightforward way to define a two-sided coideal 
$\gok_J$ in
the quantized universal enveloping algebra using the $L$-operators 
introduced in \cite{\RTF}. 
The reflection equation now more or less guarantees
the existence of sufficiently many 
$\gok_J$-fixed vectors. 
This fact is obviously central to the study of
zonal spherical functions. 
If the matrix $J$ satisfies the reflection 
equation,
the invariant ring corresponding to the coideal $\gok_J$ defines 
a quantum
homogeneous $G_q$-space that is $G_q$-isomorphic 
to the quantum orbit passing through the ``classical'' point 
associated with $J$.

In \cite{\NOmac}, the first 
author used the ideas sketched above 
to construct and analyse quantum analogues of the symmetric spaces
$GL(n)/SO(n)$ and $GL(2n)/Sp(n)$. 
The corresponding zonal spherical functions 
were shown to be expressible in terms of Macdonald's 
symmetric polynomials
corresponding to root system $A_{n-1}$ (cf.\ \cite{\MCorth}). 
As a follow-up, the first and third author extended 
this approach to all classical
compact symmetric spaces (\cite{\NS}) . 
In all these cases, 
the zonal spherical functions could
be expressed as Macdonald's polynomials (\cite{\MCorth}) or 
Koornwinder's
multivariable Askey-Wilson polynomials (\cite{\KWmaw}), 
depending on the
restricted root system of the symmetric space.

The quantum symmetric spaces treated in \cite{\NOmac}, \cite{\NS}, 
however,  do not 
exhaust all known examples of quantum analogues of classical compact
symmetric spaces. 
In the early stages of quantum group theory, Podle\'s had
already defined a continuously parametrized (parameter 
different from $q$!)
family of mutually non-isomorphic
$SU_q(2)$-homogeneous spaces. 
Each of these spaces can be regarded as a quantum
analogue of the classical 2-sphere. 
The zonal spherical functions on these
quantum spheres were analysed by Mimachi and 
the first author (\cite{\NM}),
and by \cite{\KWqsf}. They can be expressed as a certain 
subfamily of  Askey-Wilson polynomials (\cite{\AW}) or 
limit cases thereof. 
 
Recently, the first and second authors generalized these
results to quantum complex projective spaces of arbitrary dimension
(\cite{\DN}) . 
They studied a family of quantum projective spaces which were introduced
in \cite{\KV} and depend 
on a continuous parameter.
The zonal spherical functions were expressed as Askey-Wilson 
polynomials
depending on two continuous and one discrete parameter or limit cases
thereof.  See \cite{\DZ} for
a more extensive discussion of the relation between these 
parametrized families of
quantum symmetric spaces and the examples treated in 
\cite{\NOmac}, \cite{\NS}.

In this paper, we announce some results which constitute 
a generalization of \cite{\DN}. 
We present a one-parameter family of constant solutions to
one type of  reflection equation, and use them to define 
a family of quantum homogeneous $U_q(n)$-spaces that can be 
regarded as a quantum analogue
of a complex Grassmannian of arbitrary rank. 
It turns out that the zonal
spherical functions on these quantum Grassmannian spaces 
are expressed 
as a subfamily
of Koornwinder's multivariable Askey-Wilson polynomials 
depending on two continuous and two discrete parameters.
This last result is proved by computing the radial part of a 
suitable Casimir operator in the quantized universal enveloping 
algebra, a method which was first used by Koornwinder 
(\cite{\KWqsf}) in the $SU_q(2)$-case,
and subsequently generalized by the 
first 
author (\cite{\NOmac}) 
to higher rank quantum symmetric spaces.

This paper does not contain any proofs. 
Full proofs and more details will
be given elsewhere.
\vfill\eject\noindent
\beginsektion 1 {Recall on multivariable Askey-Wilson polynomials}
We recall the definition of multivariable Askey-Wilson 
polynomials (\cite{\KWmaw}). 
Our notation is  slightly
different from \cite{\KWmaw}.

Let $P_\Sigma := \bigoplus_{1\leq k \leq l} \ZZ\eps_k$
denote the weight lattice of the root system $BC_l$, 
and $P_\Sigma^+$ the cone of dominant weights in $P_\Sigma$. 
A dominant weight $\lambda = \sum_{k=1}^l \lambda_k
\eps_k\in P_\Sigma$ is then characterized by the condition
$\lambda_1 \geq \ldots \geq \lambda_l \geq 0$. Let $\leq$
denote the usual dominance ordering on weights. Recall
that $\lambda \leq \mu$ if and only if $\sum_{i=1}^k\lambda_i\leq
\sum_{i=1}^k\mu_i$ for all $1\leq k\leq l$.

We will use the notation 
$x_k := e^{\eps_k}$ ($1\le k\le\ell$)
to refer to the formal exponentials. 
Let $\CC[x^{\pm 1}]:= \CC[x_1^{\pm 1}, \ldots, x_l^{\pm 1}]$ 
denote
the algebra of Laurent polynomials in the variables 
$x_k$ ($1\leq k\leq l$). $\CC[x^{\pm 1}]$ can be viewed as the group 
algebra of the weight lattice $P_\Sigma$. 
The Weyl group $W:=\ZZ_2^l \rtimes
\goS_l$ of $BC_l$ acts on $\CC[x^{\pm 1}]$ in a 
natural way by permutations
and sign changes of the $x_k$.
 Let $\CC[x^{\pm 1}]^W\subset \CC[x^{\pm 1}]$ denote the 
subalgebra of $W$-invariant elements.

Let $0< q < 1$.
Write  $T_{q,x_k}\colon \CC[x^{\pm 1}] \to 
\CC[x^{\pm 1}]$ for the $q$-shift operator
sending $x_k$ to $qx_k$.
Koornwinder's $q$-difference operator
$D_{\eps_1}\colon \CC[x^{\pm 1}] \to \CC[x^{\pm 1}]$ 
is defined by
$$D_{\eps_1} := \sum_{k=1}^l \left ( \Phi_k^+(x) T_{q,x_k} +
\Phi_k^-(x)T^{-1}_{q,x_k} \right ) - \Phi^0(x), 
\eqno(1.1)$$
where
$$\eqalignno{\Phi_k^+(x) &:= {(1-ax_k)(1-bx_k)(1-cx_k)(1-dx_k)
\over (1-x_k^2)(1-qx_k^2)}\prod_{i\neq k} {(tx_k-x_i)
(tx_kx_i-1)\over (x_k-x_i)(x_kx_i-1)}, & \cr
\Phi_k^-(x) &:= {(x_k-a)(x_k-b)(x_k-c)(x_k-d)\over
(x_k^2-1)(x_k^2-q)} \prod_{i\neq k} {(x_k-tx_i)(x_kx_i-t)
\over (x_k-x_i)(x_kx_i-1)}, & (1.2)\cr
\Phi^0(x) &:= \sum_{k=1}^l \Phi_k^+(x) + \Phi_k^-(x). &\cr}
$$
Here we assume that the parameters $a,b,c,d,t$ are complex 
numbers satisfying
$$
0<t<1,\quad
-q\le abcd < 1.\eqno (1.3)
$$
The operator $D_{\eps_1}$ maps 
$\CC[x^{\pm 1}]^W$ into itself. 

It can be shown that there is a unique family 
$(P_\lambda)_{\lambda\in P_\Sigma^+}$ 
in $\CC[x^{\pm 1}]^W$ such that
$$
\eqalignno{\text{(i)} &\quad P_\lambda = m_\lambda + 
\sum_{\mu < \lambda}c_{\lambda\mu} m_\mu,  & (1.4)\cr
\text{(ii)} & \quad D_{\eps_1}P_\lambda = c_{\lambda\lambda}
P_\lambda. & \cr}
$$
Here $m_\lambda:= \sum_{\mu\in W \lambda} x^\mu$
is the orbit sum corresponding to $\lambda\in P^+_\Sigma$. 
One has
$$
c_{\lambda\lambda} = \sum_{k=1}^l 
\left (q^{-1} abcdt^{2l-k-1}
(q^{\lambda_k} - 1) + t^{k-1}(q^{-\lambda_k} -1)\right )
\quad (\lambda\in P_\Sigma^+).
\eqno (1.5)
$$
Note that $c_{\lambda\lambda} \neq c_{\mu\mu}$ 
for $\lambda > \mu$, provided (1.3) is 
satisfied (cf.\  \cite{\KS}, Proposition 4.6). 
The elements 
$$
P_\lambda = P_\lambda(x; a,b,c,d; q,t)
\eqno (1.6)
$$ 
will be called  {\sl multivariable 
Askey-Wilson polynomials} (associated with the root
system $BC_l$). 
For certain values of the parameters
one reobtains Macdonald's polynomials (\cite{\MCorth}) corresponding
to the pairs $(BC_l, B_l)$ or $(BC_l, C_l)$.

Assume now for simplicity 
that $0< t <1 $ and 
that 
$a,b,c,d$ are real numbers with $|a|, |b|, |c|, |d| < 1$. 
Let $T^\CC_\Sigma := (\CC^\ast)^l$ denote the complex 
torus of dimension $l$. We denote its natural compact
real form by $T_\Sigma$. Elements of $\CC[x^{\pm 1}]$ can be
naturally viewed as polynomial functions on $T^\CC_\Sigma$.

Recall the notation of {\sl $q$-shifted factorials}:
$$(a;q)_n := \prod_{k=0}^{n-1} (1-aq^k), \;
(a_1,\ldots, a_s;q)_n := \prod_{j=1}^s (a_j;q)_n,
\; (a;q)_\infty := \lim_{n\to\infty} (a;q)_n.\eqno (1.7)$$
Under the above-mentioned conditions on $a,b,c,d,t$, the
infinite product 
$$\Delta^+ := \prod_{i=1}^l {(x_i^2;q)_\infty \over
(ax_i, bx_i, cx_i, dx_i;q)_\infty} \cdot\prod_{1\leq i < j
\leq l} {(x_i/x_j;q)_\infty \over (tx_i/x_j;q)_\infty}
{(x_ix_j;q)_\infty \over (tx_ix_j;q)_\infty}\eqno (1.8)$$
defines a continuous function on $T_\Sigma$.
If we put 
$$\Delta(x):= \Delta^+(x) \Delta^+(x^{-1}),\eqno (1.9)$$
then $\Delta$ is a positive continuous function on $T_\Sigma$.
The Askey-Wilson polynomials $P_\lambda(x;a,b,c,d;q,t)$
are mutually orthogonal with respect to the weight function 
$\Delta$:
$$\int_{T_\Sigma} P_\lambda(x) P_\mu(x^{-1}) \Delta(x) \dd x
= 0, \quad \lambda\neq \mu. \eqno (1.10)$$
Here $\dd x$ denotes the Haar measure on the real torus 
$T_\Sigma$. 

For general values of the parameters, an explicit expression 
of the orthogonality measure was recently written down 
by Stokman \cite{\STaw}. 
\beginsektion 2 {Quantum complex Grassmannians}
Let us fix some notation for the quantum unitary group
(cf.\ \cite{\JB}, \cite{\RTF}, \cite{\NOmac}). 
More details can be found in
\cite{\NOmac}. 
Let $0<q<1$ and $n\geq 2$. 
Let $\Aq = \FSA_q(U(n))$  denote the
 algebra of functions on the
{\sl quantum unitary group} $U_q(n)$.
The matrix elements 
of the vector representation
$V$ with basis $(v_i)_{1\leq i\leq n}$ are written
$t_{ij}\in\Aq$ ($1\leq i,j \leq n$). They satisfy the
commutation relations 
$$R T_1 T_2 = T_2 T_1 R, \eqno (2.1)$$
where $T := (t_{ij})_{1\leq i,j \leq n}$ is an 
$n\times n$ matrix with coefficients in $\Aq$, 
$T_1 := T\ten \id$ and 
$T_2 := \id\ten T$ are Kronecker matrix products, 
and 
$R\in\End(V\ten V)$ 
 is the 
invertible $n^2\times n^2$ matrix defined by
$$R := \sum_{ij} q^{\de_{ij}} e_{ii}\ten e_{jj} +
(q-q^{-1}) \sum_{i>j} e_{ij}\ten e_{ji}.\eqno (2.2)$$
The $e_{ij}\in\End(V)$ denote the standard 
matrix units with respect to the basis $(v_i)$. We put
$$R^+ := PRP, \qquad R^- := R^{-1}, \eqno (2.3)$$
where $P\in\End(V\ten V)$ is the usual permutation operator.
The Hopf $\ast$-algebra $\Aq$
is spanned by the coefficients of its finite-dimensional unitary
corepresentations (cf.\ \cite{\WZ}, \cite{\NYM}, \cite{\DKcqg}).

Let $P$ denote the free $\ZZ$-module of rank $n$ with 
canonical basis $(\eps_i)_{1\leq i\leq n}$. 
The {\sl quantized universal enveloping algebra} $\Uq = 
\FSU_q(\gog\gol(n))$ is the algebra generated by the symbols 
$q^h$ ($h\in P^\ast$) and $e_i, f_i$ ($1\leq i \leq n-1$) 
subject to the well-known quantized Weyl-Serre relations.
The algebra $\Uq$ is also generated by  the so-called {\sl $L$-operators}
 $\Lp_{ij},\,\Lm_{ij}\in \Uq$ (cf.\ \cite{\JB}, \cite{\RTF},
\cite{\NOmac}).
The matrices $\Lpm := \sum_{ij} e_{ij} \ten \Lpm_{ij}$ with 
coefficients in $\Uq$ satisfy the following relations:
$$R^+ L_1^\epsilon L_2^\epsilon = L_2^\epsilon L_1^\epsilon R^+\; 
(\epsilon = \pm), \quad 
R^+ \Lp_1 \Lm_2 = \Lm_2 \Lp_1  R^+, \eqno (2.4)$$
where $\Lpm_1 := \Lpm \ten \id$ and 
$\Lpm_2 :=  \id \ten \Lpm$
are Kronecker matrix products.
The Hopf $\ast$-algebra structure on $\Uq$ is
determined by:
$$\De(\Lpm_{ij}) = \sum_k \Lpm_{ik} \ten \Lpm_{kj}, \quad
\eps(\Lpm_{ij}) = \de_{ij}, \quad 
(\Lpm_{ij})^\ast = S(L^\mp_{ji}),\eqno (2.5)$$
for $1\leq i,j\leq n$. We put 
$\tau = \ast \circ S\colon \Uq \to \Uq$.
 Then
$$\tau(\Lpm_{ij}) = L^\mp_{ji}\quad 
(1\leq i,j \leq n).\eqno (2.6)$$

Recall that a left $\Uq$-module $W$ is called 
{\sl $P$-weighted}
if it has a  vector space basis consisting of 
weight vectors with weights in $P$.
The cone $P^+\subset P$ of {\sl dominant weights} 
consists by definition of all weights $\lambda = 
\sum_k \lambda_k\eps_k \in P$ such that
$\lambda_1 \geq \ldots \geq \lambda_n$.
There is a 1-1 correspondence $\lambda \longleftrightarrow
V(\lambda)$ between dominant weights and irreducible
$P$-weighted finite-dimensional left $\Uq$-modules 
such that $\lambda\in P^+$
is the highest weight of $V(\lambda)$ 
(cf.\  \cite{\RS}). Recall that 
$\lambda\in P^+$ is called a highest weight of a left $\Uq$-module $W$
if there exists a non-zero vector $w\in W$ such that
$q^h\cdot w = q^{\langle h, \lambda \rangle} v$ and $e_i\cdot w = 0$
for all $1\leq i\leq n-1$. 

Given a left $\Uq$-module $W$, we define a right $\Uq$-module
structure on the same underlying vector space $W$ by 
putting
$$ 
w \cdot u := u^\ast \cdot w
\quad(w\in W, u\in\Uq).
\eqno (2.7)
$$
The corresponding right $\Uq$-module will be denoted by $W^\circ$.

The Hopf $\ast$-algebra pairing $\langle \cdot\, , 
\, \cdot \rangle$ 
between $\Uq$ and $\Aq$ is determined by 
$$\langle \Lpm_1, T_2\rangle = R^\pm, \quad 
\langle \Lpm, \detq\rangle = q^{\pm 1} \id.\eqno (2.8)$$
This pairing is nondegenerate in the sense that the 
canonical mapping 
$\Aq\to\operatorname{Hom}_{\CC}(\Uq,\CC)$ is injective. 
Under this duality, the finite-dimensional corepresentations of
$\Aq$ are in 1-1 correspondence with the finite-dimensional
$P$-weighted representations of $\Uq$.

Using the pairing $\langle \cdot\, , \, \cdot \rangle$
 between $\Uq$ and $\Aq$, one can identify $\Aq$ with a 
subspace of the algebraic linear dual of $\Uq$.
One defines
a $\Uq$-bimodule structure on $\Aq$ by putting:
$$u\cdot a := (\id \ten u) \circ \De(a), 
\quad a\cdot u := (u \ten \id) \circ \De(a)\quad 
(u\in\Uq,\, a\in\Aq).\eqno (2.9)$$
The action of the L-operators is given by:
$$L_1^\eps \cdot T_2 = T_2 R^\eps,
\quad T_2 \cdot L^\eps_1 = R^\eps T_2 \quad (\eps = \pm).
\eqno (2.10)$$
The multiplication
$\Aq \ten \Aq \to \Aq$ and the unit mapping $\CC \to \Aq$
are $\Uq$-bimodule homomorphisms. In other words,
$\Aq$ is an {\sl algebra with two-sided $\Uq$-symmetry}.
One has the following 
decomposition of $\Aq$ into irreducible $\Uq$-bimodules:
$$\Aq = \bigoplus_{\lambda\in P^+} V(\lambda) \ten V(\lambda)^\circ.
\eqno (2.11)$$
Here the subspace $W(\lambda) := V(\lambda) \ten V(\lambda)^\circ
\subset \Aq$
is spanned by the matrix coefficients of the (co-)representation
$V(\lambda)$.
The decomposition (2.11) can also be characterized 
as the simultaneous eigen\-space decomposition of $\Aq$ 
with respect to the
natural action of the center $\FSZ\Uq\subset \Uq$.
Let $h\colon \Aq\to\CC$ denote the Haar functional on $\Aq$
(cf.\ \cite{\WZ}, \cite{\DKcqg}). Then $\langle a, b\rangle :=
h(b^\ast a)$ defines a positive definite inner product on $\Aq$ 
with respect to which the subspaces $W(\lambda)\subset \Aq$ are
mutually orthogonal (Schur orthogonality). 

We now proceed to define a family of quantum Grassmannians.
Let  $X$ be an $n\times n$ matrix with coefficients in
any ring. Consider the following 
{\sl reflection
equation}:
$$R_{12}X_1 R^{-1}_{12} X_2 = X_2 R_{21}^{-1}X_1 R_{21}.
\eqno (2.12)$$
Here $X_1 := X\ten \id$ and 
$X_2 := \id\ten X$ are Kronecker matrix products,
$R_{12} = R$, $R_{21} = PRP$ ($P$ being the permutation operator).
Define the algebra $\FSC_q= \FSC_q(n)$ of functions on
the quantum space of {\sl $q$-Her\-mi\-ti\-an matrices} as the
algebra generated by the symbols $x_{ij}$ ($1\leq i,j \leq n$)
subject to the relations given by the reflection equation 
(2.12). There is a unique $\ast$-operation on $\FSC_q$ such
that $x_{ij}^\ast = x_{ji}$. In shorthand notation, we write
$X^\ast = X$.
\beginth{Proposition} two-one
{\sl There is a unique $\ast$-algebra homomorphism 
$\delta\colon \FSC_q \to \Aq \ten \FSC_q$ such that
$$\delta(x_{ij}) = \sum_{r,s} t_{ir}t^\ast_{js}\ten x_{rs}
\quad \hbox{or}\quad \delta(X)  = TXT^\ast. \eqno (2.13)$$
The mapping $\delta$ is a comodule mapping, hence defines
an action of the quantum unitary group $U_q(n)$ on the 
quantum space of $q$-Her\-mi\-ti\-an matrices.
\par}
\endth
In the terminology of \cite{\DKqsf}, a {\sl classical point}
in the quantum space of $q$-Her\-mi\-ti\-an matrices is a
$\ast$-algebra homomorphism $\tilde\eps\colon\FSC_q
\to \CC$. By the definition of the algebra $\FSC_q$, such mappings
$\tilde\eps$ are in 1-1 correspondence with $n\times n$
matrices $J$ with complex coefficients satisfying the
reflection equation (2.12) and such that $J^\ast =J$.
The correspondence is given by $\tilde\eps(X) = J$.

In the remainder of this paper, we fix an integer
$1\leq l \leq \left [ {n\over 2} \right ]$. 
Let $\sigma$ be a real parameter. 
Define an $n\times n$ matrix
$J^\sigma$ by putting
$$J^\sigma := \sum_{1\leq k\leq l} q^\sigma
(q^{-\sigma} - q^\sigma) e_{kk} + \sum_{l<k<l'} e_{kk}
-\sum_{k\leq l \ \text{or}\  k\geq l'} q^\sigma e_{kk'},
\eqno (2.14)$$
where $k':= n+1-k$ ($1\leq k\leq n$). 

\beginth{Proposition} two-two
{\sl The matrix $J^\sigma$ satisfies the reflection
equation (2.12) for any value of the parameter 
$\sigma\in \RR$.
\par}
\endth
Define a mapping
$\Psi^\sigma\colon \FSC_q \to \Aq$ by putting
$\Psi^\sigma :=  (\id \ten \tilde\eps^\sigma) \circ \delta$.
Then $\Psi^\sigma$ is a $\ast$-algebra homomorphism
intertwining the natural (left) coactions of $\Aq$ on
$\FSC_q$ and itself. We write $\FSB_q^\sigma :=
\Psi^\sigma(\FSC_q) \subset \Aq$ for the image of
$\Psi^\sigma$. $\FSB_q^\sigma$ is a $\ast$-subalgebra
and left coideal in $\Aq$. It can be viewed as the algebra
of functions on the $U_q(n)$-orbit of the classical point 
defined by $J^\sigma$ in the
quantum space of $q$-Her\-mi\-ti\-an matrices.
The action of $U_q(n)$ on this quantum orbit is transitive
by construction
(cf.\ \cite{\DKqsf}). If $q=1$ then the algebra $B_q^\sigma$
($\sigma\in \RR$)
is $U(n)$-isomorphic with the algebra of functions on the
complex Grassmannian of rank $l$ endowed with 
its natural $U(n)$-action.

The next obvious step is to look for a stabilizer ``subgroup''.
Define an $n\times n$ matrix $M^\sigma$ with coefficients in
$\Uq$ by putting
$$
M^\sigma:= \Lp J^\sigma - J^\sigma \Lm \in \End(V) \ten \Uq.
\eqno (2.15)
$$
The subspace $\gok^\sigma \subset \Uq$ is by definition
spanned by the coefficients of the matrix $M^\sigma$.
\beginth{Lemma} two-three
{\sl The subspace $\gok^\sigma\subset \Uq$ is a 
$\tau$-invariant two-sided coideal.
\par}
\endth
An element $a\in\Aq$ is called left {\sl 
$\gok^\sigma$-invariant} 
if $\gok^\sigma \cdot a = 0$.
\beginth{Proposition} two-four
{\sl The subalgebra $\FSB^\sigma_q\subset \Aq$ 
consists precisely of all left $\gok^\sigma$-invariant
elements in $\Aq$.
\par}
\endth
In other words, the two-sided coideal $\gok^\sigma$
can be regarded as an infinitesimal stabilizer of the
quantum orbit defined by the matrix $J^\sigma$.
It can be shown that, in the limit $q\to 1$, the coideal
$\gok^\sigma\subset \Uq$ tends to a Lie subalgebra of
$\gog\gol(n,\CC)$ which is conjugate to the Lie algebra
$\gok$ of the subgroup $U(l)\times U(n-l)\subset U(n)$.

Let $W$ be a finite-dimensional 
left $\Uq$-module. A vector $w\in W$ is called $\gok^\sigma$-fixed
if $\gok^\sigma\cdot w= 0$. The subspace of 
$\gok^\sigma$-fixed vectors  in $W$ is written 
$W_{\gok^\sigma}$.  We have the following key lemma:
\beginth{Lemma} two-five
{\sl Let $V$ denote the vector representation of $\Uq$,
$V^\ast$ its contragredient with dual basis $(v_i^\ast)$.
The element $w_{J^\sigma} := \sum_{i,j} J^\sigma_{ij}
v_i\ten v_j^\ast\in V\ten V^\ast$ is a $\gok^\sigma$-fixed vector.
\par}
\endth
Given any complex matrix $J$, one can define a matrix
$M_J$ with coefficients in $\Uq$ as in (2.15), and a
corresponding two-sided coideal $\gok_J\subset \Uq$.
In order that the element 
$w_J := \sum_{i,j} J_{ij}
v_i\ten v_j^\ast\in V\ten V^\ast$ ($V$ vector representation) be
a $\gok_J$-fixed vector it is necessary and sufficient for the
matrix $J$ to satisfy the reflection equation (2.12).

Let us identify elements $\lambda\in P^+$
with sequences $(\lambda_1, \ldots, \lambda_n)$ of 
 integers such that $\lambda_1 \geq \ldots
\geq \lambda_n$.
\beginth{Theorem} two-six
{\sl For any $\lambda \in P^+$, the subspace 
$V(\lambda)_{\gok^\sigma}$ is at most one-di\-men\-sion\-al.
There are non-zero $\gok^\sigma$-fixed vectors in
$V(\lambda)$ if and only if $\lambda\in P^+$ is of the form
$$\lambda = (\mu_1, \ldots, \mu_l, 0, \ldots, 0, 
-\mu_l, \ldots, -\mu_1), \quad
\mu_1 \geq \ldots \mu_l\geq 0. \eqno (2.16)$$
\par}
\endth
A dominant weight $\lambda\in P^+$ and the corresponding
representation $V(\lambda)$ are called {\sl spherical}
if $\lambda$ is of the form (2.16). The set of spherical
dominant weights will be denoted by $P_\gok^+$. 
\beginth{Corollary} two-seven
{\sl One has the following irreducible decomposition of
the right $\Uq$-module (or left $\Aq$-comodule)
$\FSB^\sigma_q$:
$$\FSB_q^\sigma = \bigoplus_{\lambda\in P^+_\gok} 
V(\lambda)^\circ. \eqno (2.17)$$
\par}
\endth
\beginsektion 3 {Zonal spherical functions}
Let $\sigma, \tau$ be real parameters.
Define $\Hst\subset \Aq$ as the subspace of left 
$\gok^\sigma$-invariant and right $\gok^\tau$-invariant 
elements. It is a $\ast$-subalgebra of $\Aq$.
\beginth{Proposition} three-one
{\sl If we put $\Hst(\lambda) := \Hst\cap W(\lambda)$
($\lambda\in P^+$) then we have the decomposition
$$\Hst = \bigoplus_{\lambda\in P^+_\gok} \Hst(\lambda).
\eqno (3.1)$$
Each of the subspaces $\Hst(\lambda)$ ($\lambda\in
P^+_\gok$) is one-dimensional.
\par}
\endth
Any non-zero element of $\Hst(\lambda)$ ($\lambda\in
P^+_\gok$) is called a {\sl zonal spherical function}
corresponding to the spherical weight $\lambda$.

Let $\AT:=\CC[z_1^{\pm 1}, \ldots, z_n^{\pm 1}]$
be the algebra of Laurent polynomials in the variables $z_i$
($1\leq i\leq n$). 
There is a unique Hopf $\ast$-algebra structure
on $\AT$ such that 
$$\Delta(z_i) = z_i \ten z_i, \quad \eps(z_i) = 1, \quad
z_i^\ast = z_i^{-1} \quad (1\leq i\leq n). \eqno (3.2)$$

There is a unique surjective Hopf $\ast$-algebra morphism
$${}_{|\TT}\colon \Aq \longrightarrow \AT$$
mapping $t_{ij}\in \Aq$  onto $\de_{ij} z_i\in \AT$. The torus
$\TT$ can be viewed as a subgroup of the quantum unitary
group $U_q(n)$. The mapping ${}_{|\TT}$ then 
is the corresponding restriction of functions.
Let us introduce new variables $x_i$ ($1\leq i\leq l$)
in the algebra $\AT$ by putting
$$x_1 := z_1z_n^{-1}, \quad x_2:= z_2z_{n-1}^{-1},
\quad \ldots,
\quad x_l := z_l z_{n+1-l}^{-1}.\eqno (3.3)$$
Note that the Weyl group $W$ of the root system
$BC_l$ acts naturally on the subalgebra of $\AT$
generated by the $x_i^{\pm 1}$. The subalgebra of $W$-invariant
elements is written 
$\CC[x^{\pm 1}]^W = \CC[x_1^{\pm 1}, \ldots, x_l^{\pm 1}]^W$.
\beginth{Theorem} three-two
{\sl The mapping ${}_{|\TT}\colon \Aq \longrightarrow \AT$
maps $\Hst$ injectively onto its image $\Hst_{|\TT}\subset
\FSA(\TT)$. One has
$$\Hst_{|\TT} = \CC[x_1^{\pm 1}, \ldots, x_l^{\pm 1}]^W.
\eqno (3.4)$$
In particular, the algebra $\Hst$ is commutative.
\par}
\endth
Consider the Casimir operator (cf. \cite{\RTF})
$$C := \sum_{ij} q^{2(n-i)}\Lp_{ij}S(\Lm_{ij})\in 
\Uq.\eqno (3.5)$$
$C$ is central in $\Uq$ and acts as a scalar on each subspace 
$W(\lambda)$ ($\lambda\in P^+$).
The corresponding eigenvalue is
$$\chi_\lambda(C) = \sum_{k=1}^n q^{2(\lambda_k + n - k)}.
\eqno (3.6)$$
Since $C$ is central in $\Uq$, the left action of
$C$ on $\Aq$ preserves the subalgebra $\Hst$, hence acts as
a scalar on each of the subspaces $\Hst(\lambda)$ 
($\lambda\in P^+_\Sigma$).
Namely, each zonal spherical function $\phi\in\Hst(\lambda) $ 
satisfies the equation 
$C\cdot\phi= \chi_\lambda(C) \,\phi$.

There is a uniquely determined linear operator
$D\colon \Hst_{|\TT} \to \Hst_{|\TT}$ 
(called {\sl the
radial part} of the Casimir operator $C$) such 
that on $\Hst$ we have 
$${}_{|\TT} \circ C = D\circ {}_{|\TT},
\eqno (3.7)$$
where the symbol $C$ denotes the 
left action of the element
$C\in\Uq$ on $\Hst\subset \Aq$.
It is possible to compute an explicit expression
for the radial part $D$. 
Recall that with any spherical weight 
$\lambda = (\mu_1, \ldots, \mu_l, 0, \ldots, 0, 
-\mu_l, \ldots, -\mu_1)$ in the weight lattice $P$
one can naturally associate
a dominant weight $\mu =  (\mu_1, \ldots, \mu_l)$
in the weight lattice $P_\Sigma$. 
\beginth{Theorem} three-three
{\sl The radial part 
$$D-\chi_\lambda(C)\id\  \colon \ \CC[x^{\pm 1}]^W \to 
\CC[x^{\pm 1}]^W
\eqno (3.8)$$
is a constant multiple of 
Koornwinder's $q$-difference operator
$D_{\eps_1} - c_{\mu\mu}\id$, defined as in  (1.1)
with base $q^2$ and parameters
$$a = -q^{\si + \tau + 1},
\; b= -q^{-\si - \tau + 1},
\; c = q^{\si - \tau + 1}, \;
d = q^{-\si + \tau + 2(n-2l) + 1}, \; t=q^2.
\eqno (3.9)$$
\par}
\endth
Note that the parameters (3.9) satisfy condition
(1.3).

For any $\lambda\in P^+_\gok$, let us a fix a non-zero
element $\phi(\lambda)\in\Hst(\lambda)$ (zonal spherical
function).
Using Theorem
(3.3) and the definition of the Askey-Wilson polynomials
$P_\mu$ we can prove:
\beginth{Theorem} three-four
{\sl The restriction 
$\phi(\lambda)_{|\TT}$ of the zonal spherical function
$\phi(\lambda)$ ($\lambda\in P^+_\gok$) to the toral
subgroup $\TT\subset U_q(n)$ is equal to the Askey-Wilson
polynomial
$$P_\mu(x; -q^{\si + \tau + 1},
-q^{-\si - \tau + 1}, q^{\si - \tau + 1}, 
q^{-\si + \tau + 2(n-2l) + 1}; q^2, q^2)\eqno (3.10)$$
up to a scalar multiple.
\par}
\endth
\beginth{Remark} three-five
As stated in section 2, the subspaces $W(\lambda)$ are
mutually orthogonal with respect to the inner product 
on $\Aq$ defined in terms of the Haar functional. This means
that the restricted zonal spherical functions 
$\phi(\lambda)_{|\TT}$ are
mutually orthogonal with respect to the inner product on
the algebra $\Hst_{|\TT}$ induced by the restriction 
of the Haar functional. 
\endth
\beginth{Remark} three-six
In the limit $\sigma \to \pm \infty$, the coideal
$\gok^\sigma$ tends to a two-sided coideal in $\Uq$
which can be viewed as a direct $q$-analogue of the
Lie subalgebra $\gog\gol(l)\oplus \gog\gol(n-l)\subset \gog\gol(n)$.
The corresponding quantum Grassmannian has the usual
Plancherel decomposition and can be regarded as the quotient
of the quantum unitary group $U_q(n)$ by the quantum 
subgroup $U_q(l) \times U_q(n-l)$. 
The limit transition of our zonal spherical functions to this case 
is expected to be consistent with the limit transition (cf.\ 
\cite{\KS})
of general multivariable Askey-Wilson polynomials to 
multivariable big $q$-Jacobi or little $q$-Jacobi polynomials
(cf.\ \cite{\STjac}). 


\Refs \widestnumber\key{NYM}
\ref\key\AW 
\by Askey, R. and Wilson, J.
\paper Some basic hypergeometric polynomials 
that generalize Jacobi polynomials
\jour Mem. Amer. Math. Soc. \vol 54 \yr 1985 \issue 319
\endref

\ref\key\DZ 
\by Dijkhuizen, M.S. \yr 1994
\paper Some remarks on the construction
of quantum symmetric spaces 
\paperinfo 
Proceedings of a conference on Representation Theory of Lie groups,
Lie algebras, and their Quantum Analogues, Twente (NL), 
Dec. 1994 
\jour Acta Appl. Math. \toappear\endref

\ref\key \DKqsf 
\manyby 
Dijkhuizen, M.S., and Koornwinder, T.H.
\paper Quantum homogeneous spaces, duality 
and quantum 2-spheres
\jour Geom. Dedicata \vol 52 \yr 1994 \pages 
291-315\endref 

\ref\key \DKcqg 
\bysame 
\paper CQG algebras: a direct algebraic approach
to compact quantum groups \jour Lett. Math. Phys. \vol 32
\yr 1994 \pages 315-330\endref

\ref\key\DN
\by Dijkhuizen, M.S., Noumi, M. 
\paper A family of quantum projective spaces 
and related $q$-hyper\-geo\-metric 
orthogonal polynomials \paperinfo preprint \yr 1995\endref

\ref\key \DR
\by Drinfel'd, V.G. \paper Quantum groups 
\paperinfo in:  Proceedings ICM Berkeley (1986), ed. A.M. Gleason
\publ Amer. Math. Soc. \publaddr Providence, RI \yr 1986
\pages 798-820\endref

\ref\key \JB 
\by Jimbo, M. 
\paper A $q$-analogue of $U(\gog\gol(n))$,
Hecke algebra and the Yang-Baxter equation 
\jour Lett. Math. Phys.
\vol 11 \yr 1986 \pages 247-252\endref

\ref\key\KWqsf 
\manyby Koornwinder, T.H.
\paper Askey-Wilson polynomials as zonal
spherical functions on the $SU(2)$ quantum group
\jour SIAM J. Math. Anal. \vol 24 \issue 3 \yr 1993
\pages 795-813 \endref

\ref\key\KWmaw 
\bysame \paper Askey-Wilson
polynomials for root systems of type $BC$
\inbook in:\ ``Hypergeometric functions on domains of
positivity, Jack polynomials, and applications'',
ed. D.S.P. Richards, Contemp. Math. 138
\publ Amer. Math. Soc. \publaddr Providence, RI
\yr 1992\pages 189-204\endref

\ref\key\KV
\by Korogodsky, L.I., and Vaksman L.L.
\paper Quantum $G$-spaces and Heisenberg algebra 
\inbook in:\hfill\break ``Quan\-tum Groups'',
ed. P.P. Kulish, Lecture Notes in Math. 1510 
\publ Springer-Verlag \yr 1992 \pages 56-66\endref

\ref\key\KU
\by Kulish, P.P. \paper Quantum groups and quantum
algebras as symmetries of dynamical systems
\paperinfo preprint YITP/K-959 \yr1991 \endref

\ref\key\MCorth 
\by Macdonald, I.G. \paper Orthogonal
polynomials associated with root systems \paperinfo
preprint \yr 1988\endref

\ref\key \NOmac 
\by Noumi, M. \paper 
Macdonald's symmetric polynomials
as zonal spherical functions on some 
quantum homogeneous spaces
\jour Adv. Math. \toappear \yr 1993\endref

\ref\key \NM 
\by Noumi, M., and Mimachi, K. \paper
Quantum 2-spheres and big $q$-Jacobi polynomials
\jour Comm. Math. Phys. \vol 128 \yr 1990 \pages 521-531
\endref

\ref\key \NS 
\by Noumi, M., and Sugitani, T.  \yr 1995
\paper Quantum symmetric spaces and related 
$q$-orthogonal polynomials \inbook in:\ 
``Group Theoretical Methods in Physics'' \bookinfo 
Proceedings XX ICGTMP, Toyonaka (Japan), 1994,
ed. A. Arima et al. \publ World Scientific
\publaddr Singapore \pages 28-40\endref

\ref\key \NYM
\by Noumi, M., Yamada, H., and Mimachi, K. 
\paper Finite dimensional representations of the quantum 
group $GL_q(n;\Bbb{C})$ and the zonal spherical functions 
on $U_q(n-1)\backslash U_q(n)$
\jour Japanese J. Math \vol 19 \yr 1993
\pages 31-80\endref

\ref\key\PO 
\by Podle\'s, P. \paper Quantum spheres 
\jour Lett. Math. Phys. \vol 14 \yr 1987 
\pages 193-202\endref

\ref\key\RTF 
\by Reshetikhin, N., Faddeev, L.D., and Takhtajan, L.A. 
\paper Quantization of Lie groups 
and Lie algebras \jour Leningrad
Math. J. \vol 1 \yr 1990 \pages 193-225\endref

\ref\key\RS 
\by Rosso, M. 
\paper Finite-dimensional representations of
the quantum analog of a complex simple Lie algebra 
\jour Comm. Math. Phys.
\vol 117 \yr 1988 \pages 581-593\endref

\ref\key\STjac 
\manyby Stokman, J.V.
\paper Multivariable big and little $q$-Jacobi polynomials
\paperinfo Mathematical Preprint Series 95-16,
University of Amsterdam \yr 1995 \endref

\ref\key\STaw
\bysame \paper Multivariable BC type Askey-Wilson polynomials
with partly discrete orthogonality measure \paperinfo
preprint 
\yr 1995
\endref

\ref\key\SK
\by Stokman, J.V.,  and Koornwinder, T.H.
\paper Limit transitions for BC type multivariable 
orthogonal polynomials \paperinfo Mathematical Preprint
Series 95-19, University of Amsterdam \yr 1995 \endref

\ref\key\WZ 
\by Woronowicz, S.L. \paper Compact matrix pseudogroups
\jour Comm. Math. Phys. \vol 111 \yr 1987 
\pages 613-665\endref

\endRefs
\enddocument